\documentstyle[12pt]{amsart}
\makeatletter
\def\LaTeX{\leavevmode L\raise.42ex    \hbox{\kern-.3em\size{\sf@size}{0pt}
\selectfont A}\kern-.15em\TeX}
\makeatother

\newcommand{\BibTeX}{{\rm B\kern-.05em{\sc i\kern-.025emb}\kern-.08em\TeX}}

\input xy
\xyoption{all}

\newtheorem{thm}{Theorem}[section]
\newtheorem{lem}[thm]{Lemma}
\newtheorem{prop}[thm]{Proposition}
\newtheorem{cor}[thm]{Corollary}
\theoremstyle{definition}
\newtheorem{defn}{Definition}[section]
\newtheorem{exmp}{Example}[section]
\makeatletter\label{e:dispaa}
\label{e:dispau}
\label{e:dispav}
\label{e:dispaw}
\label{e:dispax}
\def\theequation{\thesection.\@arabic\c@equation}
\makeatother
\makeatletter

\newcommand{\lra}{\longrightarrow}

\newcommand{\benu}{\begin{enumerate}}
\newcommand{\enu}{\end{enumerate}}

\makeatother
\begin{document}
\title[postprojectives and components]%
{On the postprojective partitions and components of the Auslander-Reiten quivers}

\author[Coelho]
{Fl\'avio U. Coelho}
\address{Departamento de Matem\'atica-IME, 
Universidade de S\~ao Paulo, CP 66281, 
S\~ao Paulo, SP, 05315-970, Brazil}
\email{fucoelho@@ime.usp.br}

\author[Silva]
{Danilo D. da Silva}
\address{Departamento de Matemática-DMA, 
Universidade Federal de Sergipe,
S\~ao Crist\'ov\~ao, SE, 49100-000, Brazil}
\email{ddsilva@@ufs.br}
\maketitle

\begin{abstract}
In this paper we shall investigate further the connections between the 
postprojective partition of an algebra and its Auslander-Reiten quiver. 
\end{abstract}
\vspace{.3 cm}

Auslander-Smal\o\ introduced, in \cite{auslander}, the notion of postprojective partition and modules
(under the name of preprojective). The connection between such a partition and the structure of the 
Auslander-Reiten quiver has been investigated in several papers such as \cite{assemcoelho,auslander,
coelho1, coelho2,coelhosilva1,todorov}. The purpose of this paper is to follow such investigations.

We introduce the notion of ${\bf P}$-discrete component of the Auslander-Reiten quiver $\Gamma_A$ as 
follows. Let $\{ \bf{P_i} \}$ of ${\rm ind}A$ with $i\in \mathbb{N}_{\infty}=\mathbb{N}\cup{\{ \infty} \}$
be the  postprojective partition of $A$ (recall the definition below). A component $\Gamma$ of $\Gamma_A$ 
is ${\bf P}$-discrete if for each $i\geq 0$ and each $M \in \Gamma \cap {\bf P_i}$, we have that 
tr$_{\bf P_{i+1}}(M) = tr_{\bf P_\infty} (M)$, where tr$_{\mathcal C}(M)$ denotes the trace of the set of 
modules ${\mathcal C}$ in $M$ (see Section 1 below). 
\vspace{.3 cm}\\
{\bf Theorem 2.3.} 
{\it
Let $A$ be a representation-infinite Artin algebra. If $\Gamma$ is a $\bf{P}$-discrete connected component of $\Gamma_A$
then there is no arrow $ M \rightarrow N$ in $\Gamma$ with $M \in \bf{P_i}$ and $N \in \bf{P_j}$ such that $i+1<j<{\infty}$.
}

\vspace{.3 cm}
Also, using the notion of left degree of a morphism (introduced by Liu in \cite{liu}), we prove 
the following result.
\vspace{.3 cm}\\
{\bf Theorem 3.4.} 
{\it
Let $A$ be a finite dimensional algebra over an algebraically closed field and let 
$f \colon M \longrightarrow N$ be an irreducible monomorphism of infinite left degree. If tr$_{\bf P_\infty}(N) = 0$, 
then ${\rm Coker}f\in{\bf P}_{\infty}$ and every non-trivial submodule of $ {\rm Coker}f$ is postprojective. 
}
\vspace{.3 cm}

In particular, from the last theorem we get that if an irreducible monomorphism which lies in a postprojective component has its cokernel in a regular component then the latter must be simple regular.  

This paper is organized as follows. After recalling basic notions in Section 1, we prove Theorem 2.3 in 
Section 2 and Theorem 3.4 in Section 3.

\section{Preliminaries}

\subsection{Basics}
For the results of Section 2, we assume  the algebras $A$ to be Artin algebras, unless otherwise stated. For the last section, we shall restrict to finite dimensional algebras over a fixed algebraically closed field $k$. Futhermore, we will assume that all algebras are basic.

For unexplained notions in representation theory we refer the reader to \cite{auslanderbook}.

For an algebra $A$, we denote by $\rm{mod}A$ the category of all finitely generated left $A$-modules, and by $\rm{ind}A$ the full subcategory of $\rm{mod}A$ consisting of one representative of each isomorphism class of indecomposable $A$-modules. We denote by $\Gamma_A$ the Auslander-Reiten quiver of $A$ and by $\tau$ and $\tau^{-}$ the Auslander-Reiten translations DTr and TrD, respectively.

Given $n \geq 1$ and  $M,N \in {\rm mod}A$, we define the subgroups ${\rm rad}^n(M,N)$ of ${\rm Hom}(M,N)$ by induction: for $n=1$, we set ${\rm rad}^1(M,N)$ to be the set of all morphisms $f \colon M \lra N$ such that the compositions $gfh$ are not 
isomorphisms for all $h \colon L \lra M$ and $g \colon N \lra L$, with $L$ indecomposable. Also, we define ${\rm rad}^n(M,N)$ as the set of all morphisms $f \in {\rm Hom}(M,N)$ such that there exist $X \in {\rm mod}A$ and morphisms $g \in {\rm rad}(M,X)$ and $h \in {\rm rad}^{n-1}(X,N)$ such that $f=hg$. Finally, we set ${\rm rad}^{\infty}(M,N)= \bigcap_{n \geq 1}{\rm rad}^n(M,N)$.

We recall that for $X,Y \in {\rm ind}A$, $f:X \rightarrow Y$ is called {\bf irreducible} if and only if $f \in {\rm rad}(X,Y) \backslash {\rm rad}^2(X,Y)$. A {\bf path of irreducible morphisms of length n} is a sequence $ M_0 \stackrel{h_1}{\lra} M_1 \lra \cdots \lra M_{n-1} \stackrel{h_n}{\lra} M_n$ where each $h_i$ is irreducible and each $M_j$ is indecomposable.

Following Liu \cite{liu}, we say that the left degree of an irreducible morphism $f:X \rightarrow Y$ is $n$, and we denote $d_l(f)=n$, if $n$ is the smallest positive integer for which there exist $Z \in {\rm ind}A$ and a morphism $h : Z \rightarrow X$ such that  $h \in {\rm rad}^n(Z,X) \backslash{\rm rad}^{n+1}(Z,X)$ and $fh \in {\rm rad}^{n+2}(Z,Y)$. In case this condition is not verified for any $n \geq 1$  we say the left degree of $f$ is infinite. Dually, one can define the right degree of an irreducible morphism.

\subsection{Postprojective partitions and modules} 
We shall now recall the concept of postprojective partition and modules as introduced by Auslander and Smal\o\ in \cite{auslander} under the name preprojective. 

A {\bf postprojective partition} of an Artin algebra $A$ is a partition $\{ \bf{P_i} \}$ of ${\rm ind}A$ with $i\in
\mathbb{N}_{\infty}=\mathbb{N}\cup{\{ \infty} \}$ such that

\begin{itemize}
\item[(a)] ${\rm ind}A$ is the disjoint union of the subcategories ${\bf
P_i}$, $i \in \mathbb{N}_{\infty}$.
\item[(b)] for each $j<{\infty}$, ${\bf P_j}$ is a finite minimal cover of the union of the subcategories ${\bf
P_i}$ such that $j \leq i \leq {\infty}$. 
\end{itemize}

It is clear that the $A$-modules in ${\bf P_0}$ are all the indecomposable projectives. In this article, we denote ${\bf P}({\rm ind}A)= \displaystyle\bigcup_{0 \leq i <{\infty}}{\bf P_i}$ simply by ${\bf P}$. The modules in ${\rm add}{\bf P}$ will be called {\bf postprojective} modules (former preprojective modules in \cite{auslander}). We denote by ${\bf P^m}$ the subcategory ${\bf P_0} \cup \cdots \cup {\bf P_{m}}$.

Given $i \in \mathbb{N}_{\infty}$ and a module $M$ in ${\rm mod}A$
we denote the trace of ${\bf P_i}$ on $M$ by ${\rm tr}_{\bf P_i}(M)$, that is, the submodule of $M$ generated by the images of all morphisms which have domain in ${\rm add}{\bf P_i}$. Therefore, ${\rm tr}_{\bf P_i}(M)$ is the submodule of  $M$ generated by $\{ {\rm Im}f | f \in {\rm Hom}(N,M) \, {\rm and} \, N \in {\bf P_i} \}$. It was proved in \cite{auslander} that ${\rm tr}_{\bf P_{\infty}}(M)= \displaystyle\cap_{i\geq 0}
{\rm tr}_{\bf P_{i}}(M)$. Hence, $ {\rm tr}_{\bf P_{\infty}}(M)= {\rm tr}_{\bf 
P_{r}}(M)$ for some $r \in \mathbb{N}$ since $M$ is artinian and ${\rm tr}_{\bf P_{n+1}}(M) \subseteq {\rm tr}_{\bf P_{n}}(M)$,  for each $n \geq 0$. It was also proved in \cite{auslander} that $M$ is postprojective if and only if ${\rm tr}_{\bf P_{\infty}}(M) \neq M$.

The following proposition from \cite{coelho1} shall be very useful in the sequel.

\begin{prop}
\label{lema1}
{\rm \cite{coelho1}}
Let $N$ be a postprojective module and $f : M \rightarrow N$ a morphism such that $ {\rm Im}f \not\subseteq {\rm tr}_{\bf{P_{\infty}}}(N)$. Then $f \not\in {\rm rad}^{\infty}(M,N)$.
\end{prop}

We also recall the following result from \cite{coelhosilva1} (Lemma 4.2).

\begin{lem}
\label{lemadan}
Let  $ \bf{P_0},\bf{P_1}, \cdots, {\bf P_{\infty}}$ be the postprojective partition of an algebra $A$. Given $0<i \leq {\infty}$, we have ${\rm Hom}(M,N)={\rm rad}^i(M,N)$, for each $ M \in \bf{P_0}$ and each $ N \in \bf{P_i}$.
\end{lem}

\section{Path of irreducible morphisms and the postprojective partition}

Along this section, let $A$ denote an Artin algebra and $\{ \bf{P_i} \}$ of ${\rm ind}A$ with $i\in
\mathbb{N}_{\infty}=\mathbb{N}\cup{\{ \infty} \}$ the postprojective partition of ind$A$. Let $M$ be 
a postprojective module. It was proved  in \cite{auslander} that there exists a path of irreducible morphisms 
from a projective module $P$ to $M$. Clearly, then, $P$ and $M$ lie in the same connected component 
$\Gamma$ of $\Gamma_A$. One could wonder if there exists such a path as follows:
$$ P = M_0 \lra M_1 \lra \cdots \lra M_n= M$$
with $M_i \in {\bf P_i}$ for each $i$. Corollary 3 in \cite{igusa} states that this is true if $\Gamma$ is 
a postprojective component of a hereditary algebra. The next example show that this is 
not true in general. , and, on the other hand, that there are non-postprojective components with this 
property. 

\begin{exmp}
{\rm
Let $A$ be the finite-dimensional $k$-algebra (where $k$ is a field) given by the quiver\\
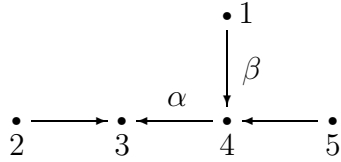
\begin{picture}(50,20)
\put(43,0){\circle*{1}}
\put(57,0){\circle*{1}}
\put(71,0){\circle*{1}}
\put(85,0){\circle*{1}}
\put(71,14){\circle*{1}}

\put(45,0){\vector(1,0){10}}
\put(69,0){\vector(-1,0){10}}
\put(83,0){\vector(-1,0){10}}
\put(71,12){\vector(0,-1){10}}

\put(63,2){$\alpha$}
\put(73,5.6){$\beta$}

\put(72.5,13){1}
\put(42,-4.5){2}
\put(56,-4.5){3}
\put(70,-4.5){4}
\put(84,-4.5){5}

\end{picture}
\vspace{1 cm}\\
bound by $\beta \alpha = 0$. Its Auslander-Reiten quiver has the following shape:\\
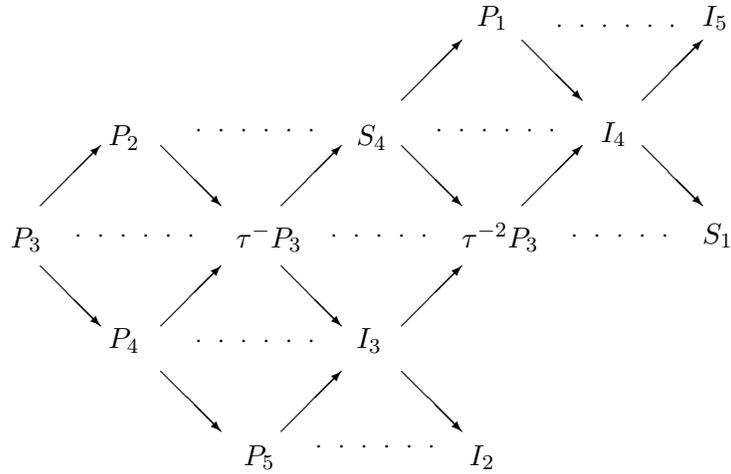
\begin{picture}(50,40)
\put(52,4){\vector(1,1){8}}
\put(68,12){\vector(1,-1){8}}
\put(52,-4){\vector(1,-1){8}}
\put(68,-12){\vector(1,1){8}}

\put(62,12){{\small $S_4$}}
\put(46,-1.5){{\small $\tau^{-}P_3$}}
\put(62,-15){{\small $I_3$}}
\put(76,-1.5){{\small $\tau^{-2}P_3$}}

\put(20,4){\vector(1,1){8}}
\put(84,4){\vector(1,1){8}}
\put(68,18){\vector(1,1){8}}
\put(100,18){\vector(1,1){8}}
\put(36,-12){\vector(1,1){8}}
\put(52,-26){\vector(1,1){8}}

\put(36,12){\vector(1,-1){8}}
\put(100,12){\vector(1,-1){8}}
\put(84,26){\vector(1,-1){8}}
\put(20,-4){\vector(1,-1){8}}
\put(36,-18){\vector(1,-1){8}}
\put(68,-18){\vector(1,-1){8}}

\put(16,-1.5){{\small $P_3$}}
\put(108,-1){{\small $S_1$}}
\put(29,-15){{\small $P_4$}}
\put(47,-30.5){{\small $P_5$}}
\put(77,-30.5){{\small $I_2$}}
\put(29,12){{\small $P_2$}}
\put(94.5,12.5){{\small $I_4$}}
\put(78,28){{\small $P_1$}}
\put(108,28){{\small $I_5$}}

\multiput(59,0)(3,0){5}{\circle*{.1}}
\multiput(25,0)(3,0){6}{\circle*{.1}}
\multiput(91,0)(3,0){5}{\circle*{.1}}
\multiput(41,14)(3,0){6}{\circle*{.1}}
\multiput(73,14)(3,0){6}{\circle*{.1}}
\multiput(41,-14)(3,0){6}{\circle*{.1}}
\multiput(57,-28)(3,0){6}{\circle*{.1}}
\multiput(89,28)(3,0){6}{\circle*{.1}}

\end{picture}
\vspace{4 cm}\\
For each $j$, $P_j, I_j$ and $S_j$ denote, respectively, the projective, the injective and 
the simple modules associated to the vertex $j$ of the quiver. 
The postprojective partition is then ${\bf P_0} = \{ P_1, P_2, P_3, P_4, P_5 \}$, 
${\bf P_1} = \{ \tau^{-1} P_3, I_3, I_4 \}$, ${\bf P_2} = \{ S_4, \tau^{-2} P_3, I_1, I_2 \}$ and 
${\bf P_3} = \{ I_5 \}$. Observe that there are paths from a projective to $I_5 \in {\bf P_3}$ of 
length 2, 4 , 5 and 6, and so, none of the required type. 
}
\end{exmp}

Next example shows that there are non-postprojective components with this property. 

\begin{exmp}
\label{tuboraio} 
{\rm 
Let $A$ be a path algebra defined by the quiver:

$$
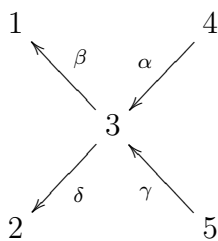
\xymatrix{
       &          1      &                              &          4     \ar[ld]_{\alpha}                                  \\
         &                &                  3     \ar[ld]^{\delta} \ar[lu]_{\beta}      &                                                   \\
         &       2         &                              &          5 \ar[lu]^{\gamma}            \\
}$$

\par\noindent
(see \cite{coelhosilva1}). 
Let $S$ be the simple module associated to the vertex $3$ and let $N$ be the indecomposable module such that $\tau S=N$ and $\tau N= S$. One can check that $S$ and $N$ determine a tube of rank 2.\\
The bounded quiver of the extended algebra $A[S]$ is

$$
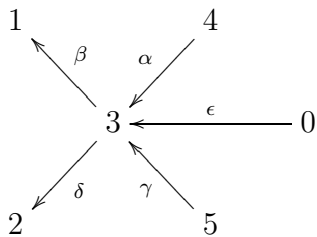
\xymatrix{
       &          1      &                              &          4     \ar[ld]_{\alpha}           &                       \\
         &                &                  3     \ar[ld]^{\delta} \ar[lu]_{\beta}      &      &   0 \ar[ll]_{\epsilon}                                          \\
        &       2         &                              &          5 \ar[lu]^{\gamma}      &      \\
}$$

\par\noindent
bounded by $\epsilon\beta=0$ and $\epsilon\delta=0$.\\
The ray tube which contains $S=S[1]$ has the shape:

$$
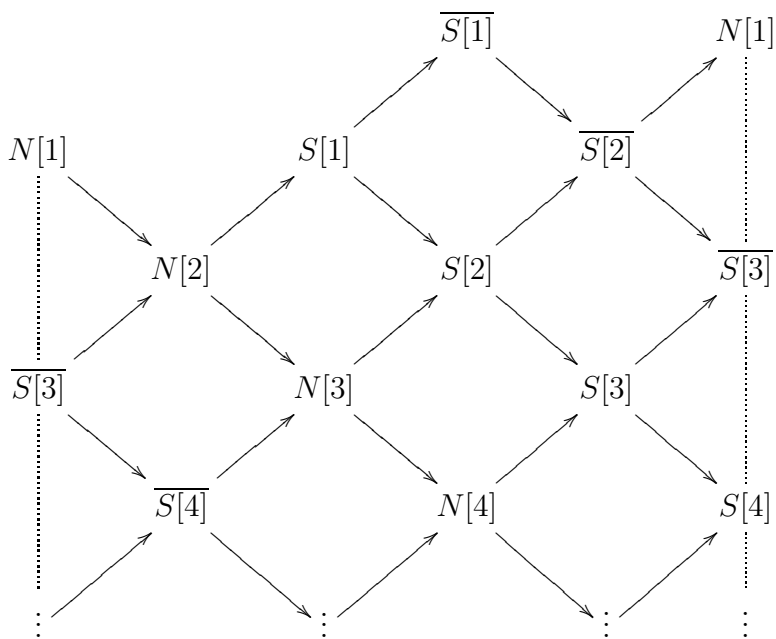
\xymatrix{
                                                                &                                                      &                                                                   &   \overline{S[1]}  \ar[rd]                    &                                                       &        N[1]   \ar@{.}[dd]                      \\
N[1]  \ar[rd]  \ar@{.}[dd]                               &                                                      & S[1]  \ar[rd] \ar[ru]                                   &                                                            &  \overline{S[2]} \ar[ru] \ar[rd]      &                                                   \\
                                                                &   N[2]   \ar[rd] \ar[ru]                  &                                                                   &  S[2]   \ar[rd] \ar[ru]                          &                                                       &    \overline{S[3]} \ar@{.}[dd]             \\
\overline{S[3]} \ar[rd] \ar[ru] \ar@{.}[dd]    &                                                      &    N[3]    \ar[rd] \ar[ru]                             &                                                           &   S[3] \ar[rd] \ar[ru]                       &                                           \\
                                                                &  \overline{S[4]} \ar[rd] \ar[ru]     &                                                                  &  N[4]   \ar[rd] \ar[ru]                         &                                                        &   S[4]      \ar@{.}[d]                         \\
         \vdots       \ar[ru]                                     &                                                      &     \vdots   \ar[ru]                                               &                                                            &  \vdots \ar[ru]                                         &              \vdots                            \\
}$$

\vspace{0.3cm}

\par\noindent
Let ${\bf P_0}, {\bf P_1}, \cdots, {\bf P_{\infty}}$ be the postprojective partition of $A[S]$. Then the modules $N[i]$ and $S[j]$ with $i \geq 1$ and $j \geq 1$ are all in ${\bf P_{\infty}}$ and the ray $\overline{S[1]} \rightarrow \overline{S[2]} \rightarrow \cdots \rightarrow \overline{S[n]} \rightarrow \cdots$ is such that $\overline{S[n]} \in {\bf P_{n-1}}$ for each $n \geq 1$. 
}
\end{exmp}

Our purpose in this section is to introduce the concept of $\bf{P}$-discrete components of $\Gamma_A$ and show that, 
for them, we have an affirmative answer to the above question. We start with a lemma.

\begin{lem}
\label{formadopreprojetivo}
Let $n$ be an integer greater than 0, and let $M\in {\bf P_n}$. For each $j$, $0 \leq j < n$, there exists a path of irreducible morphisms 
$L \leadsto M$ through postprojective modules with $L \in {\bf P_j}$ and with composition not lying in ${\rm rad}^{\infty}(L,M)$.
\end{lem}
\begin{pf} 
Since $M\in{\bf P_n}$ and $j <n$, by definition, there exists an epimorphism 
$\displaystyle\bigoplus_{i=1}^{r}L_i \stackrel {g}{\longrightarrow} M$, with $L_i \in {\bf P_j}$ for each $i$. Write $g=[g_1, \cdots, g_r]$. 
Since $M$ is postprojective, tr$_{\bf P_{\infty}}(M) \neq M$ and so Im$g \not\subseteq {\rm tr}_{\bf P_{\infty}}(M)$. Hence, there exists 
an $l$ such that Im$g_l \not\subseteq {\rm tr}_{\bf P_{\infty}}(M)$. Because of Proposition \ref{lema1}, $g_l \not\in {\rm rad}^{\infty} (M_l, M)$. 
By \cite{auslanderbook} (Proposition 7.4), we get the expression $g_l = 
\sum_i \alpha_i + \beta $ where each $\alpha_i$ is a path of irreducible morphisms with composite not 
lying in rad$^{\infty} (L_l, M)$ and $\beta \in$ rad$^{\infty} (L_l, M)$. Using Proposition \ref{lema1} again, we infer that 
Im$\beta \subseteq$tr$_{\bf P_{\infty}}(M)$. Hence, for some $i$, the composition $h$ of the path $\alpha_i$ has image 
Im$h\not\subseteq {\rm tr}_{\bf P_{\infty}}(M)$ and so $h \notin$rad$^{\infty} (L_l, M)$. It remains to show that the path $\alpha_i$ pass 
through only postprojective modules. Suppose $\alpha_i$ is a path 
$$ L_l \stackrel{(*)}{\leadsto} N \stackrel{(**)}{\leadsto} M$$
where $N \in {\bf P_{\infty}}$ and write by $\gamma$ and $\gamma'$  the (nonzero) compositions of the paths $(*)$ and $(**)$, respectively. 
Hence $\gamma ' \gamma = h$. Now, Im$\gamma' \subseteq {\rm tr}_{\bf P_{\infty}}(M)$ because $N\in {\bf P_{\infty}}$ and so 
Im$ h \subseteq$ Im$\gamma' \subseteq {\rm tr}_{\bf P_{\infty}}(M)$, a contradiction. This proves the lemma. 
\end{pf}

\begin{defn}
Suppose  $A$  is representation-infinite. We say that a connected component $\Gamma$ of $\Gamma_A$ is a $\bf{P}$-{\bf{discrete component}} if for all $i \geq 0$ and for each postprojective module $M$ in $\Gamma$ with $M \in \bf{P_i}$ we have ${\rm tr}_{\bf P_{i+1}}(M)={\rm tr}_{\bf P_{\infty}}(M)$.
\end{defn}

\par\noindent
{\bf Remark:} Note that ${\rm tr}_{\bf P_{i+1}}(M)={\rm tr}_{\bf P_{\infty}}(M)$ implies ${\rm tr}_{\bf P_{j}}(M)={\rm tr}_{\bf P_{\infty}}(M)$, for each $j>i$. 

\begin{prop}
\label{caracterizacao}
Let $A$ be a representation-infinite Artin algebra and $\Gamma$ be a connected component of $\Gamma_A$. The following are equivalent:

\begin{enumerate}
\item[(a)] $\Gamma$ is a $\bf{P}$-discrete component.
\item[(b)] There exists no arrow $M \rightarrow N$ in $\Gamma$ with $M \in \bf{P_j}$, $N \in \bf{P_i}$ and $i<j<{\infty}$.
\item[(c)] There exists no path of irreducible morphisms $M \leadsto N$ in $\Gamma$ through indecomposable postprojective modules with $M \in \bf{P_j}$, $N \in \bf{P_i}$ and $i<j<{\infty}$.
\end{enumerate}
\end{prop}

\begin{pf}
As one can easily see that (b) and (c) are equivalent we show (a)$\Rightarrow $(b)  and (c)$\Rightarrow $(a).

(a)$\Rightarrow $(b) 
Suppose there exists an irreducible morphism $f:M \rightarrow N$ in $\Gamma$ with $M \in \bf{P_j}$, $N \in \bf{P_i}$ and $i<j<{\infty}$. Assuming  $\Gamma$ is $\bf{P}$-discrete we have ${\rm Im}f \subseteq {\rm tr}_{\bf P_{j}}(N) = {\rm tr}_{\bf P_{\infty}}(N)$. Then we can factorize $f$ as follows:
$$\xymatrix{
M \ar[rr]^f  \ar[rd]_{f'}  & & N \\
     & {\rm tr}_{\bf P_{\infty}}(N) \ar @{^{(}->} [ru] &  \\
}$$
From the fact that $f$ is irreducible and $N \neq {\rm tr}_{\bf P_{\infty}}(N)$, we get that $f':M \rightarrow {\rm tr}_{\bf P_{\infty}}(N)$ is a split monomorphism and $M$ is a summand of ${\rm tr}_{\bf P_{\infty}}(N)$. This is an absurd since ${\rm tr}_{\bf P_{\infty}}(N) \in {\rm add}\bf{P_{\infty}}$.

(c)$\Rightarrow $(a)  
Suppose by contradiction that there exists $M \in \Gamma \cap \bf{P_i}$, with ${\rm tr}_{\bf P_{i+1}}(M) \neq {\rm tr}_{\bf P_{\infty}}(M)$. Then there exists $f: M' \rightarrow M$, $M' \in \bf{P_{i+1}}$, such that ${\rm Im}f \not\subseteq {\rm tr}_{\bf P_{\infty}}(M)$. Therefore, Lemma \ref{formadopreprojetivo} provides us a path of irreducible morphisms through indecomposable postprojective modules starting at $M'$ and ending at $M$ which contradicts (c). 
\end{pf}


\begin{thm}
\label{theotodorov}
Let $A$ be a representation-infinite Artin algebra. If $\Gamma$ is a $\bf{P}$-discrete connected component of $\Gamma_A$
then there is no arrow $ M \rightarrow N$ in $\Gamma$ with $M \in \bf{P_i}$ and $N \in \bf{P_j}$ such that $i+1<j<{\infty}$.
\end{thm}
\begin{pf}
We shall prove it by induction on $i \geq 0$. \\
For $i = 0$, just observe that if $f: P \rightarrow N$ is an irreducible morphism in $\Gamma$ with $P \in \bf{P_0}$ and $N \in \bf{P_j}$, $j>1$, then by Lemma \ref{lemadan} we have $f \in {\rm rad}^2(P,N)$ which is a contradiction.\\
Suppose now that the theorem is true for all values less than $i$ and let  $f:M \rightarrow N$ be an irreducible morphism in $\Gamma$ with $M \in \bf{P_i}$, $N \in \bf{P_j}$ and $i+1<j<{\infty}$. Then $\tau N \in {\bf P^{i-1}}$ (Lemma 2.1 in \cite{coelho1}) and $f$ is not a sink morphism because otherwise, by Lemma \ref{formadopreprojetivo}, there would be a path of irreducible morphisms through indecomposable postprojective modules $L \rightarrow \cdots \rightarrow M \rightarrow N$, with $L \in \bf{P_{i+1}}$ which contradicts the fact that $\Gamma$ is ${\bf P}$-discrete. Therefore, there exists $M' \in \ $mod$A$ such that

$$
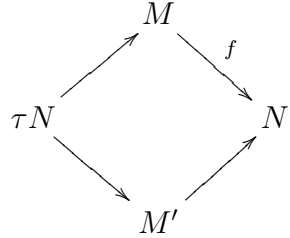
\xymatrix{
            & M \ar[rd]^{f} &  \\
  \tau N \ar[rd] \ar[ru] &    & N \\ 
            & M' \ar[ru] &    \\
}$$
is the Auslander-Reiten sequence which ends at $N$.\\
We now prove that $M' \not\in \ $add$\bf{P^{j-2}}$. Suppose, by contradiction, that $M' \in \ $add$\bf{P^{j-2}}$ and hence $M \oplus M' \in \ $add$\bf{P^{j-2}}$. Then there exists $h: N' \rightarrow N$, $N' \in \bf{P_{j-1}}$, such that ${\rm Im}h \not\subseteq {\rm tr}_{\bf P_{\infty}}(N)$, $h \in {\rm rad}(N',N)$. Since $h$ is not a split epimorphism, it can be lifted through $f$: 
$$\xymatrix{
  & N' \ar[d]^h \ar[ld]_g & \\
M \oplus M' \ar[r]_f & N \ar[r]  &  0 \\
}$$
Since $\Gamma$ is ${\bf P}$-discrete, we know that ${\rm tr}_{\bf P_{j-1}}(M \oplus M')={\rm tr}_{\bf P_{\infty}}(M \oplus M')$ which implies 
that ${\rm Im}g \subseteq {\rm tr}_{\bf P_{\infty}}(M \oplus M')$. Therefore ${\rm Im}fg \subseteq {\rm tr}_{\bf P_{\infty}}(N)$ since $f({\rm tr}_{\bf P_{\infty}}(M \oplus M')) \subseteq {\rm tr}_{\bf P_{\infty}}(N)$. Hence ${\rm Im}h = {\rm Im}fg \subseteq {\rm tr}_{\bf P_{\infty}}(N)$, which is a contradiction. \\
Then there exist a summand $M_1 \in \bf{P_{j-1}}$ of $M'$ and an irreducible morphism $\tau N \rightarrow M_1$ such that $ \tau N \in \bf{P^{i-1}}$ and $i-1<i<j-1$ which contradicts the induction hypothesis. The theorem follows.
\end{pf}

\par\noindent
{\bf Remark:} In \cite{todorov}( item (a) of Proposition 1), it was proved that if $A$ is a hereditary algebra then there is no arrow $ M \rightarrow N$ in the postprojective component with $M \in \bf{P_i}$ and $N \in \bf{P_j}$ such that $i+1<j<{\infty}$. Then this result was used to prove that the postprojective component satisfies item (b) of Proposition \ref{caracterizacao} (in other words, the converse of Theorem \ref{theotodorov} for $A$ hereditary).

\begin{cor}
Let $A$ be a representation-infinite Artin algebra and $\Gamma$ be a $\bf{P}$-discrete component of $\Gamma_A$. Given $i>0$ and $M \in \bf{P_i} \cap \Gamma$, then there exists a path of irreducible morphisms between indecomposable modules  $M_0 \lra M_1 \lra \cdots \lra M_i = M $, where $M_j \in {\bf P_j}$ for each $j \in \{0, \cdots, i \}$, and moreover it has the smallest possible length among all the paths of irreducible morphisms starting at a projective and ending at $M$.
\end{cor}

\begin{pf} 
Consider the sink map ending at $M$. Since $M$ is not projective, this morphism is an epimorphism. Then there exists an irreducible morphism $N \lra M$ with $N \in \bf{P^{i-1}}$. In fact, by Theorem 2.3, we do have that $N \in \bf{P_{i-1}}$. We keep using the same argument until we get a path of irreducible morphisms $M_0 \lra M_1 \lra \cdots \lra M_i = M $, where $M_j \in {\bf P_j}$ for each $j \in \{0, \cdots, i\}$. Again by Theorem 2.3 we get that there can not be a smallest path starting at a projective and ending at $M$.
\end{pf}

\section{Right degrees of irreducible morphisms in $\pi$-components}

Recall that if an irreducible morphism $f:M \lra N$ in ${\rm mod}A$ has finite right degree then $f$ is a monomorphism and we have $d_r(f)=n$ if, and only if, ${\rm coker}(f) \in {\rm rad}^n \backslash {\rm rad}^{n+1}$(by the dual version of Corollary 3.3 in \cite{chaio}). We are now particularly interested in a particular case, assuming in addition that ${\rm tr}_{\bf P_{\infty}}(N)=0$. We look for a connection between the fact that the $A$-module ${\rm Coker}f$ is postprojective to the fact that the right degree of $f$ is finite. Recall that, in \cite{coelhosilva1}, we have proved that the inclusion map $f_S: {\rm rad}P_S \hookrightarrow P_S$, where $P_S$ is the projective covering of the simple $S$,  has finite right degree if and only if $S$ is postprojective. From now on, since we depend on the results of \cite{chaio}, we shall restrict our consideration to finite dimensional algebras over an algebraically closed field $k$. Unless otherwise stated, $A$ is such an algebra. 

\begin{thm}
\label{grauengana}
Let $f:M \rightarrow N$ be an irreducible monomorphism with ${\rm tr}_{{\bf P_{\infty}}}(N)=0$. Then $d_r(f)<{\infty}$ if and only if ${\rm Coker}f $ is  postprojective. \label{propmono}
\end{thm}

\begin{pf}
Suppose ${\rm Coker}f$ is postprojective. Then by the fact that ${\rm coker}(f)$ is an epimorphism and ${\rm tr}_{\bf P_{\infty}}({\rm Coker}f) \neq {\rm Coker}f$ we have ${\rm coker}(f) \not\in {\rm rad}^{\infty}(N,{\rm Coker}f)$, by Proposition \ref{lema1}, so we get $d_r(f)<{\infty}$.\\
Now assume $d_r(f)=n$, $1 \leq n < {\infty}$, and suppose by contradiction $C ={\rm Coker}f \in {\bf P_{\infty}}$.\\ 
We set $\pi_f={\rm coker}(f)$. Then $\pi_f \in {\rm rad}^n(N,C ) \backslash {\rm rad}^{n+1}(N,C)$, by Proposition 3.5 in \cite{chaio}. We know there exists $1 \leq j < {\infty}$ such that ${\rm tr}_{\bf P_{j}}(N) = {\rm tr}_{\bf P_{\infty}}(N)=0$. Moreover, there exists a nonzero morphism $v:L \rightarrow C$ with $L \in {\rm add}{\bf P_j}$ such that $v \in {\rm rad}^{n+1}(L,C)$. Indeed, if we take $r>j+n$ then we can get a covering $h_r:M_r \rightarrow C$ with $h_r \in {\rm rad}(M_r, C)$ and $M_r \in {\rm add}{\bf P_r}$, since $C \in {\bf P_{\infty}}$. Then let $h_l:M_l \rightarrow M_{l+1}$ be a covering with $M_l \in {\rm add}{\bf P_l}$ for all $l \in \{ j, \cdots ,r-1 \}$. If we set $v=h_rh_{r-1} \cdots h_j:M_j \rightarrow C$ then we have that $v$ is nonzero as a composition of epimorphisms and $v \in {\rm rad}^{n+1}(L,C)$ with $L=M_j \in {\bf P_j}$ as required.  By Proposition 5.6 in \cite{auslanderbook}, either there exists $p:N \rightarrow L$ such that $\pi_f=vp$ or there exists $q:L \rightarrow N$ such that $v=\pi_fq$. In the first case, we get $\pi_f \in {\rm rad}^{n+1}(N,C)$ as $v \in {\rm rad}^{n+1}(L,C)$, which is a contradiction. In the last case, as ${\rm tr}_{\bf P_{k}}(N) = {\rm tr}_{\bf P_{\infty}}(N)=0$ we have $q=0$ which implies $v=0$, again a contradiction. Hence ${\rm Coker}f$ is postprojective. 
\end{pf}

In \cite{coelho2}, Coelho has considered the so-called $\pi$-components in $\Gamma_A$, which are components containing only postprojective modules. Hereditary algebras, or more generally, left glued algebras (see \cite{assemcoelho}), contain such components. These components can also be characterized for the fact that all its modules $M$ satisfies ${\rm tr}_{\bf P_{\infty}}(M)=0$ (see \cite{coelho2}). 

\begin{cor}
Let $\Gamma$ be a $\pi$-component and  $f:M \rightarrow N$ be an irreducible monomorphism in $\Gamma$. Then $d_r(f)<{\infty}$ if and only if ${\rm Coker}f $ is  postprojective.
\end{cor} 

\par\noindent
{\bf Remark:} Although the above results may suggest that it is true that if $f:M \rightarrow N$ is an irreducible monomorphism with $N$ postprojective then $d_r(f)<{\infty}$ if and only if ${\rm Coker}f \in {\bf P}$, we alert that it is not the case. In Example \ref{tuboraio}(b) we see a source map $f:\overline{S[1]} \rightarrow \overline{S[2]}$ which is a monomorphism of right degree equal to 1. We have that the $A[S]$-module $\overline{S[2]}$ is postprojective but $N[1]={\rm Coker}f \in {\bf P}_{\infty}$.
\vspace{ .3 cm}

Now we turn our attention to the case $d_r(f)={\infty}$. We start with a definition. 

\begin{defn}
Let $M$ be an indecomposable in ${\bf P_{\infty}}$. We say $M$ is ${\bf P_{\infty}-}{\rm simple}$ if every nontrivial submodule of $M$ is postprojective.
\end{defn}

\par\noindent
{\bf Remark:} We know that there is no bound on the lengths of the modules lying in any given infinite set of 
indecomposable nonisomorphic postprojective modules (see \cite{auslander}). Hence, any $A$-module $M$ must have a finite number of 
nonisomorphic postprojective submodules (if any). Therefore if $M$ is ${\bf P_{\infty}-}{\rm simple}$ then there exists $0 \leq n < {\infty}$ such that every nontrivial submodule of $M$ is in ${\rm add}{\bf P^n}$.

\begin{prop}
Let $M$ be and indecomposable $A$-module in a regular component of $\Gamma_A$. If $M$ is ${\bf P}_{\infty}$-simple then $M$ is simple regular.
\end{prop}
\begin{pf}
If M is ${\bf P}_{\infty}$-simple then $M$ must not have nontrivial regular submodules since all modules in regular components are in ${\bf P}_{\infty}$ . 
\end{pf}

Now we show that if $f:M \rightarrow N$ is an irreducible monomorphism with ${\rm tr}_{{\bf P_{\infty}}}(N)=0$ such that $d_r(f)={\infty}$ then ${\rm Coker}f$ is ${\bf P_{\infty}-}{\rm simple}$. By Theorem in \cite{coelho2}, we know that every monomorphism in a $\pi$-component satisfies this condition.

\begin{thm}
\label{simpleregular}
Let $f:M \rightarrow N$ be an irreducible monomorphism with ${\rm tr}_{{\bf P_{\infty}}}(N)=0$. If $d_r(f)={\infty}$ then ${\rm Coker}f$ is ${\bf P_{\infty}-}{\rm simple}$.
\end{thm}

\begin{pf}
We know by Theorem \ref{propmono} that ${\rm Coker}f \in {\bf P_{\infty}}$. Suppose by contradiction that ${\rm Coker}f$ is not ${\bf P_{\infty}-}{\rm simple}$. Then there exists a nontrivial submodule $X$ of ${\rm Coker}f$ such that $X \in {\bf P_{\infty}}$. Hence ${\rm Hom}(X,N)=0$ since ${\rm tr}_{{\bf P_{\infty}}}(N)=0$. Consider the exact sequence $0 \rightarrow M \stackrel{f}{\rightarrow} N \stackrel{g}{\rightarrow} {\rm Coker}f \rightarrow 0$ and the inclusion morphism $v:X \hookrightarrow {\rm Coker}f$. Then by Proposition 5.6 in \cite{auslanderbook} either there exists $q:X \rightarrow N$ such that $v=gq$ or there exists $p:N \rightarrow X$ such that $g=vp$. The former leads to a contradiction because ${\rm Hom}(X,N)=0$ implies $q=0$ which in turn implies $v=0$ and $X=0$. The latter also leads to a contradiction for $g$ being an epimorphism implies that $v$ is an epimorphism and hence $X={\rm Coker}f$ which is not possible. Therefore ${\rm Coker}f$ is ${\bf P_{\infty}-}{\rm simple}$.

\end{pf}

\begin{cor}
\label{corregular}
If $f:M \rightarrow N$ is an irreducible monomorphism in a $\pi-component$ such that ${\rm Coker}f$ lies in a regular component then ${\rm Coker}f$ is simple regular.\\
\end{cor}

We give now an example which illustrates the above corollary.

\begin{exmp}{\rm Let $A$ be a path algebra defined by the quiver:

$$
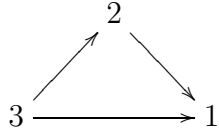
\xymatrix{
                       &                2      \ar[rd]   &                                \\
      3 \ar[rr]  \ar[ru]    &                            &  1                            \\
}$$

For each vertex $x$ in the quiver we set $P_x$ the correspondent projective. Then the irreducible monomorphism $f:P_1 \rightarrow P_3$ has infinite right degree since ${\rm Coker}f$ is not in the postprojective component. ${\rm Coker}f$ is as follows:

$$
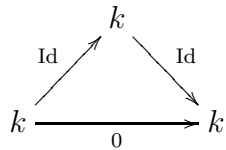
\xymatrix{
                       &                k      \ar[rd]^{\rm Id}   &                                \\
      k \ar[rr]_0  \ar[ru]^{\rm Id}    &                            &  k                            \\
}$$

On the other hand, let $\mu \in {\rm rad}^2(P_1, P_3)$ be the composition of the irreducibles $P_1 \rightarrow P_2 \rightarrow P_3$. Then $f'=f+ \mu$ is also irreducible and ${\rm Coker}f'$ is as follows:

$$\xymatrix{
                       &                k      \ar[rd]^{\rm Id}   &                                \\
      k \ar[rr]_{\rm Id}  \ar[ru]^{\rm Id}    &                            &  k                            \\
}$$

Since $f$ has infinite right degree we have that $f'$ has also infinite right degree. One can easily verify that ${\rm Coker}f$ and ${\rm Coker}f'$ lie in two different homogeneous tubes. By the above corollary both modules should be simple regular modules. And in fact that is what happens since both ${\rm Coker}f$ and ${\rm Coker}f'$ lie in the mouth of its respectives tubes.}

\end{exmp}

\subsection*{ACKNOWLEDGEMENTS} This paper is part of the PhD's thesis of the second author \cite{tese}, written under the supervision of the 
first. They both acknowledge financial suppport from CNPq during the preparation of this work. 


\end{document}